\documentclass[12pt]{article}
\usepackage{amsmath,amssymb,amsbsy,amsfonts,amsthm,latexsym,
            amsopn,amstext,amsxtra,euscript,amscd,amsthm}

\newtheorem{lem}{Lemma}
\newtheorem{lemma}[lem]{Lemma}

\newtheorem{prop}{Proposition}
\newtheorem{proposition}[prop]{Proposition}

\newtheorem{thm}{Theorem}

\newtheorem{cor}{Corollary}
\newtheorem{corollary}[cor]{Corollary}

\newtheorem{defi}{Definition}
\newtheorem{definition}[defi]{Definition}

\newtheorem{rem}{Remark}

\def\\{\cr}
\def\({\left(}
\def\){\right)}
\def\[{\left[}
\def\]{\right]}
\def\<{\langle}
\def\>{\rangle}

\begin{document}

\title{On Permutation Binomials over Finite Fields}

\author{{\sc Mohamed~Ayad} \\
{Laboratoire de Math{\'e}matiques Pures et Appliqu{\'e}es}\\
{Universit\'e du Littoral}\\
{F-62228 Calais, France} \\
{ayad@lmpa.univ-littoral.fr}
\and
{\sc Kacem~Belghaba} \\
{D\'epartment de Math\'ematiques}\\
{Universit\'e d'Oran \`a Es Senia}\\
{Bp 1524, Algeria} \\
{belghaba.kacem@univ-oran.dz}
\and
{\sc Omar~Kihel} \\
{Department of Mathematics}\\
{Brock University} \\
{Ontario, Canada L2S 3A1} \\
{okihel@brocku.ca}}

\date{\today}

\pagenumbering{arabic}

\maketitle
\begin{abstract}
Let $\mathbb{F}_{q}$ be the finite field of characteristic $p$ containing $q = p^{r}$ elements and $f(x)=ax^{n} + x^{m}$ a binomial with coefficients in this field. If some conditions on the gcd of $n-m$ an $q-1$ are satisfied then this polynomial does not permute the elements of the field. We prove in particular that if $f(x) = ax^{n} + x^{m}$ permutes $\mathbb{F}_{p}$, where $n>m>0$ and $a \in {\mathbb{F}_{p}}^{*}$, then $p -1  \leq (d -1)d$, where $d = {\mbox{gcd}}(n-m,p-1)$, and that this bound of $p$ in term of $d$ only, is sharp. We show as well how to obtain in certain cases a permutation binomial over a subfield of $\mathbb{F}_{q}$ from a permutation binomial over $\mathbb{F}_{q}$.
\end{abstract}


\section{Introduction}
Let $\mathbb{F}_{q}$ be the finite field of characteristic $p$ containing $q = p^{r}$ elements. A polynomial $f(x) \in \mathbb{F}_{q}$ is called a permutation polynomial of $\mathbb{F}_{q}$ if the induced map $f: \mathbb{F}_{q} \rightarrow \mathbb{F}_{q}$ is one to one. The study of permutation polynomials goes back to Hermite \cite{Her} for $\mathbb{F}_{p}$ and Dickson \cite{Dic} for $\mathbb{F}_{q}$. The interest on permutation polynomials increased in part because of their application in cryptography and coding theory. Despite the interest of numerous people on the subject, characterizing permutation polynomials and finding new families of permutation polynomials remain open questions. Carlitz conjectured that given an even positive integer $n$ there exists a constant $C(n)$ such that for $q > C(n)$, then there are no permutation polynomials of degree $n$ over $\mathbb{F}_{q}$. Fried, Guralnick and Saxl \cite{Fri} proved Carlitz's conjecture. Permutation monomials are completely understood, however permutation binomials are not well understood. Niederreiter and Robinsom \cite{Nie} proved the following theorem.

\begin{thm}
Given a positive integer $n$, there is a constant $C(n)$ such that for $q > C(n)$, no polynomial of the form $ax^{n} + bx^{m} + c \in \mathbb{F}_{q}[x]$, with $n > m> 1$, $\gcd(n, m) = 1$, and $ab \neq 0$, permutes $\mathbb{F}_{q}$.
\label{thm:Neiderreiter}
\end{thm}

The constant $C(n)$ in Theorem 1 is not explicit. Turnwald \cite{Tur} improved Theorem 1 and proved the following.
\begin{thm}
If $f(x) = ax^{n} + x^{m}$ permutes $\mathbb{F}_{q}$, where $n>m>0$ and $a \in {\mathbb{F}_{q}}^{*}$, then either $q \leq (n-2)^{4} + 4n -4$ or $n = mp^{i}$. 
\end{thm}
Turnwald's proof uses Weil's lower bound \cite{Wei} for the number of the points on the curve $(f(x)-f(y)/(x-y)$ over $\mathbb{F}_{q}$. For $q$ a prime number, Turnwald \cite{Tur} proved the following.
\begin{thm}
If $f(x) = ax^{n} + x^{m}$ permutes $\mathbb{F}_{p}$, where $n>m>0$ and $a \in {\mathbb{F}_{p}}^{*}$, then $p < n{\mbox{max}}(m,n-m)$. 
\end{thm}
For $m =1$, Wan \cite{Wan2} proved the following.
\begin{thm}
If $f(x) = ax^{n} + x$ permutes $\mathbb{F}_{p}$, where $n>1$ and $a \in {\mathbb{F}_{p}}^{*}$, then  $p -1 \leq (n-1) \cdot{\mbox{gcd}}(n-1,p-1)$.
\end{thm}
The bounds in Theorem 3 and Theorem 4 are of different nature. The bound in Theorem 3 is given in term of ${\mbox{max}}(m,n-m)$, whereas the bound in Theorem 4 is given in term of ${\mbox{gcd}}(n-1,p-1)$.
Theorem 3 and Theorem 4 have been improved (\cite{Mas1}) as follow. 
\begin{thm}
If $f(x) = ax^{n} + x^{m}$ permutes $\mathbb{F}_{p}$, where $n>m>0$ and $a \in {\mathbb{F}_{p}}^{*}$, and $d = {\mbox{gcd}}(n-m,p-1)$, then $p -1  \leq (n -1) \cdot {\mbox{max}}(m,d)$. 
\end{thm}
The bounds in the theorems above are not given in term of $d$ only, and one can ask whether the prime $p$ can be bounded in term of $d$ only. The answer was given by Masuda and Zieve \cite{Mas1} who proved the following.
\begin{thm}
If $f(x) = ax^{n} + x^{m}$ permutes $\mathbb{F}_{p}$, where $n>m>0$ and $a \in {\mathbb{F}_{p}}^{*}$. Then $p -1  \leq (d +1)d$. 
\end{thm}
Clearly, Theorem 6 improves Theorem 5 whenever $d-1 \leq n-1$, which is always the case except when $m =1$ and $(n-1) \mid (p-1)$. In section 2, we prove the following.
\begin{thm}
If $f(x) = ax^{n} + x^{m}$ permutes $\mathbb{F}_{p}$, where $n>m>0$ and $a \in {\mathbb{F}_{p}}^{*}$. Then $p -1  \leq (d -1)d$. 
\end{thm}
Clearly, Theorem 7 implies Theorem 6 and Theorem 5 in all cases. When $m =1$ and $(n-1) \mid (p-1)$, we will see in corollary 5 that $ p-1 \leq (n-1)(n-3)$, which improve Theorem 5. 
It would be interesting to have a bound for $p$ in term of $d = {\mbox{gcd}}(n-m,p-1)$ when $f(x) = ax^{n} + x^{m}$ permutes $\mathbb{F}_{q}$ and $q$ a power of the prime $p$. In Theorem 9, we will show how in certain cases, one can obtain from a permutation binomial $f(x) \in \mathbb{F}_{q}[x]$ a new permutation binomial $g(x) \in \mathbb{F}_{p}[x]$, and deduce in Corollary 6,  a bound of $p$ in term of $d = {\mbox{gcd}}(n-m,p-1)$. Some consequences of this theorem are stated in sections $2$ and $3$

We fix some notation which will be used through this paper. The letter $p$ always denotes a prime number, and $\mathbb{F}_{q}$ the finite field containing $q = p^{r}$ elements. For any polynomial $g(x) \in \mathbb{F}_{q}[x]$, we denote by $\overline{g(x)}$ the unique polynomial of degree at most $q-1$, with coefficients in $\mathbb{F}_{q}$ such that $g(x) \equiv \overline{g(x)} \pmod{(x^{q}-x)}$. When we refer  to a binomial $f(x)$ over $\mathbb{F}_{q}$, we always mean a polynomial $f(x) \in \mathbb{F}_{q}[x]$ of the form $f(x) = ax^{n} + x^{m}$ with  the nonrestrictive condition $\gcd(m, n) = 1$ (see \cite[Ex. 2. 1]{Sma1}), $n>m$, and $a \neq 0$. The integer $d = \gcd(n-m, q-1)$ will play an important role. It is well known that if $-a\in(\mathbb{F}_{q}^{\star})^d$ , then the equation $f(x)=0$ has $d+1$ distinct solutions in $\mathbb{F}_{q}$, hence $f(x)$ is not a permutation of $\mathbb{F}_{q}$\cite{Sma2}. In particular this claim is true if $d=1$. 

\section{Non existence of Permutation Binomials of Certain Shapes}

An old and stricking result in the theory of permutation polynomials, is the following theorem proved by Hermite for the prime fields and Dickson in the general case. 

\begin{thm}
\label{thm:Hermite/Dickson}
Let $p$ be a prime number, $q = p^{r}$, and $g(x) \in \mathbb{F}_{q}[x]$. Then $g(x)$ is a permutation polynomial if and only if
\begin{itemize}
\item [(i)] $g(x) = 0$ has a unique solution in $\mathbb{F}_{q}$.
\item[(ii)] For every $l \in \lbrace 1, \ldots, q-2 \rbrace, {\mbox{deg}}\,\overline{g^{l}(x)} \leq q-2$.
\end{itemize}
\end{thm}

For binomials, we deduce from Theorem \ref{thm:Hermite/Dickson} the following corollary.
\begin{cor}
Let $f(x) = ax^{n} + x^{m} \in \mathbb{F}_{q}[x]$, such that $a\neq 0$ and $(m, n) = 1$. Let $d = \gcd(n-m, q-1)$. Suppose $d \geq 2$. Then $f(x)$ is a permutation polynomial of $\mathbb{F}_{q}$ if and only if
\begin{itemize}
\item [(i)] $f(x) = 0$ has a unique solution in $\mathbb{F}_{q}$.
\item[(ii)] For every $l \in \lbrace 1, \ldots, q-2\rbrace$ such that $d \mid l$, we have deg $\overline{f^{l}(x)} \leq q - 2$.
\end{itemize}
\end{cor}

\begin{proof}
From Theorem \ref{thm:Hermite/Dickson}, we have only to prove that if $l \in \lbrace 1,\ldots, q-2 \rbrace$ and $d \nmid l$, then deg $\overline{f^{l}(x)} \leq q-2$. Let $k$ be an integer and let $\overline{k}$ be the integer in $\lbrace1,\ldots,q-1 \rbrace$ such that $k \equiv \overline{k} \pmod{q-1}$. Then

\begin{equation*}
 x^{k} \equiv \left\lbrace  \begin{array}{ll}
         1 & \mbox{if } k = 0\\
        x^{\overline{k}} & \mbox{if } k \neq 0 \end{array}\right.  
\end{equation*}

It follows that if $k > 0$, then $x^{k} \equiv x ^{q-1} \pmod{x^{q}-x}$ if and only if $k \equiv 0 \pmod{q-1}$. Suppose that there exists $l \in \lbrace 1,\ldots, q-2 \rbrace$ with $d \nmid l$ such that deg $\overline{f^{l}(x)} = q-1$. We deduce from 
\begin{align}
(ax^{n}+x^{m})^l & = \sum^{l}_{j=0}\binom{l}{j}a^{j}x^{nj+m(l-j)}\nonumber\\
& = \sum^{l}_{j=0}\binom{l}{j}a^{j}x^{(n-m)j+lm}
\label{eq: 1}
\end{align}
 that there exists an integer $j \in \lbrace 0, \ldots, l \rbrace$ such that
 \begin{equation*}
 x^{(n-m)j + lm} \equiv x^{q-1} \pmod{x^{q}-x}
 \end{equation*}
 Hence, $(n-m)j + lm > 0$ and $(n-m)j+lm \equiv 0 \pmod{q-1}$. Since $d = \gcd(n-m, q-1)$, then $d \mid (n-m)$ and $d \mid q-1$. But $gcd(n, m) = 1$ implies that $\gcd(d, m) = 1$. Then $d \mid l$ which is a contradiction.
\end{proof}

Corollary 1 reduces enormously the calculations when applying Theorem 8 to check whether a given a binomial  permutes  $\mathbb{F}_{q}$ or not. One needs to check the degress of only $ \[\frac{q-2}{d}\]$ polynomials instead of $q-2$ polynomials as given by Theorem 8 (ii).

For the proof of Theorem 7, we need the following lemma.

\begin{lem}
Let $f(x)$ be a binomial such that $d > 1$. Let $l \in \lbrace 1, \ldots, q-2 \rbrace$ such that $d \mid l$. Then the following assertions are equivalent

\begin{itemize}
\item[(i)] \begin{equation} {\mbox{deg}}f^{l}(x) \leq q-2 \label{eq:2} \end{equation}
\item[(ii)] \begin{equation} \sum_{\substack{j=0\\(n-m)j + lm\equiv 0\pmod{q-1}}}^{l}\binom{l}{j}a^{j} = 0 \label{eq:3} \end{equation}
\item[(iii)] \begin{equation} \sum_{\lambda = 0}^{\gamma_{l}} \binom{l}{j_{0}+\lambda (q-1)/d} \left( a^{(q-1)/d} \right)^{\lambda}  = 0 \label{eq:4} \end{equation}
 where $j_{0}$ is the smallest integer $\geq 0$ satisfying 
 \begin{equation*} j_{0} \equiv \frac{- lm}{(n-m)} \pmod{\frac{q-1}{d}} \equiv \frac{-lm}{d} / \frac{n-m}{d} \pmod{\frac{q-1}{d}} \end{equation*} 
 and $\gamma_{l}$ is the largest  integer $\lambda$ such that \begin{equation*} j_{0} + \lambda(q-1)/d \leq l \end{equation*}
\end{itemize}
\end{lem}

\begin{proof}
From equation 1, deg $\overline{f^{l}(x)} \leq q-2$ if and only if
\begin{equation}
\sum_{\substack{j=0\\(n-m)j + lm\equiv 0\pmod{q-1}}}^{l}\binom{l}{j}a^{j} = 0
\label{eq:5}
\end{equation}
The condition $(n-m) j + lm \equiv 0 \pmod{q-1}$ is equivalent to $\frac{(n-m)}{d}j + \frac{l}{d}m \equiv 0 \pmod{\frac{q-1}{d}}$, which is equivalent to
\begin{equation}
j \equiv \frac{-lm}{(n-m)} \pmod{\frac{q-1}{d}}
\label{eq: 6}
\end{equation}
Let $j_{0}$ be the smallest integer satisfying (6). Then $j \equiv j_{0} \pmod{(q-1)/d}$. Hence, equation (5) is equivalent to 
\begin{equation*}
\sum_{\lambda = 0}^{\alpha_{l}} \binom{l}{j_{0} + \lambda (\frac{a-1}{d})} \left( \left( a \right)^{\frac{q-1}{d}} \right)^{\lambda} = 0 
\end{equation*}
where $\gamma_{l}$  is the largest integer $\lambda$ such that $j_{0} + \lambda \left( \frac{q-1}{d}\right) \leq l$.
\end{proof}

\noindent {\bf Proof of Theorem 7.} In the following proof, we repetedly use Equation (4). The integer $j_{0}$ appearing in this equation depends on $l$. So, it will be denoted by $j_{0}(l)$. Suppose that there exists a permutation  binomial $f(x)$ over $\mathbb{F}_{p}$ such that $p - 1 >  d(d-1)$, then $(p-1)/d > d-1$, i.e., $(p-1)/d \geq d$ . From Theorem 10, $(p-1)/d  \neq d$, then  $(p-1)/d > d$. Set $(p-1)/d = \alpha d - z$ where $\alpha > 1$ and $z$ are integers such that $z \in \lbrace 0, \ldots, d-1 \rbrace$. By Theorem 10, we may suppose that $z \in \lbrace 1, \ldots, d-1 \rbrace$. Let $j_{0}(d)$ be the unique integer determined by (4) for $l =  d$. Set $j_{0}(d) = \beta d + \delta$ with $\delta \in \lbrace 0, \ldots, d-1 \rbrace$, then
\begin{equation}
j_{0}(d) < (p-1)/d < \alpha d.
\end{equation}

\begin{enumerate}
\item[(Case 1):] $j_{0}(d) \leq d$.

In this case, because $j_{0}(d) + \lambda(p-1)/d > d$ for $\lambda \geq 1$,  Equation (4) reduces to $\binom{d}{j_{0}(d)} = 0$. Since $j_{0}(d) \leq d < p$, then $\binom{d}{j_{0}(d)} \neq 0$, which is a contradiction, and we can exclude this case.

\item[(Case 2):] $j_{0}(d) > d$.

Clearly, $\beta \geq 1$, and From (7), we deduce that $\beta < \alpha$. Consider  Equation (4) for $l = \alpha d$. We have 
\begin{align*}
\alpha j_{0}(d)& = (\frac{p-1}{d^{2}} + \frac{z}{d})j_{0}(d)\\
& = \frac{p-1}{d}\beta + z\beta + \alpha\delta \\
&\equiv z\beta + \alpha\delta \pmod{(p-1)/d}.
\end{align*}

\item[(Case 2.1):] $z\beta + \alpha\delta < (p-1)/d$.

In this case, we have
\begin{equation}
d <  j_{0}(\alpha d) = z \beta + \alpha\delta < (p-1)/d < \alpha d = l.
\end{equation}
Let $\lambda$ be a positive integer, then 
\begin{align*}
j_{0}(\alpha d) + \lambda(p-1)/d & \geq j_{0}(\alpha d) + (p-1)/d  \\
&= z \beta + \alpha\delta + (p-1)/d \\
&> d + (p-1)/d \\
&> z + (p-1)/d \\
&= \alpha d 
= l.
\end{align*}
Hence there is only one term in the left hand side of equation (4) corresponding to $l = \alpha d$, namely $\binom{\alpha d}{j_{0}(\alpha d)} = \binom{\alpha d}{z\beta + \alpha \delta}$.  Since $(p-1)/d \geq d$, we have $ \alpha d < \frac{p-1}{d} + d < p $. Hence, from (8), we obtain that $j_{0}(\alpha d) = z \beta + \alpha\delta  < \alpha d < p$. Then $\binom{\alpha d}{j_{0}(\alpha d)} = \binom{\alpha d}{z\beta + \alpha \delta} \neq 0$,  and we reject this case.
 
\item[(Case 2.2):] $z\beta + \alpha\delta \geq (p-1)/d$.

Suppose that $\delta = 0$, then $z\beta + \alpha\delta = z\beta$ and since $\beta \leq \alpha - 1$, we deduce that $z\beta + \alpha\delta \leq (\alpha - 1)z \leq(\alpha - 1)(d-1)$, hence 
\begin{align*}
z\beta + \alpha\delta &\leq\alpha d - \alpha - d + 1\\
& < \alpha d - z = (p-1)/d,
\end{align*}
which is a contradiction. We many suppose that $\delta$ is positive. Consider  Equation (4) for  $l = (\alpha - 1)d$. We have 
\begin{align*}(\alpha - 1)j_{0}(d) &= (\frac{p-1}{d^{2}} + \frac{z-d}{d})j_{0}(d)\\
& = (\frac{p-1}{d^{2}} + \frac{z-d}{d})(\beta d + \delta) \\
&\equiv (z-d)\beta + (\alpha - 1)\delta \pmod{(p-1)/d}.
\end{align*}

In order to prove that  $(z-d)\beta + (\alpha-1)\delta =j_{0}\left( (\alpha - 1)d \right)$, we have to show that  $0 \leq (z-d)\beta + (\alpha-1)\delta < (p-1)/d$. Since $z < d$, then $(z-d)\beta < 0$, hence 
\begin{align*}
(z-d)\beta + (\alpha - 1)\delta &< (\alpha - 1)\delta \leq(\alpha - 1)(d-1) \\
&= \alpha d - d - \alpha + 1 \\
&< \alpha d - z  = (p-1)/d.
\end{align*}

We now look at the sign of $(z-d)\beta + (\alpha - 1)\delta$.

On the one hand side, we have $\alpha > \beta \geq 1$, hence $\alpha \geq 2$. Furthermore since $z - d < 0$ and $\beta \leq \alpha - 1$, then $(z-d)\beta \geq (z-d)(\alpha - 1)$, hence 
\begin{align*}
(z-d)\beta + (\alpha - 1)\delta &\geq (z-d)(\alpha - 1) + (\alpha - 1)\delta \\
&= (\alpha - 1)(z-d+\delta) \\
&\geq z- d + \delta.
\end {align*}

On the other hand side, since $z\beta + \alpha\delta \geq (p-1)/d = \alpha d - z$, then 
\begin{align*}
(z-d)\beta + (\alpha - 1)\delta &= z\beta + \alpha\delta - d\beta - \delta \\
&\geq \alpha d - z - d\beta - \delta \\
&= (\alpha  -\beta)d  - z - \delta \\
&\geq d  - z - \delta.
\end {align*}
We have  shown that $(z-d)\beta + (\alpha - 1)\delta \geq |A|$, where $A = z - d - \delta$. Hence  $(z - d)\beta + (\alpha - 1)\delta \geq 0$ and then $j_{0}\left( (\alpha - 1)d \right) = (z-d)\beta + (\alpha - 1)\delta$. As in the preceding cases we prove that in the left hand side of equation (4), for $l = (\alpha - 1)d$, there is only one term. For any integer $\lambda \geq 1$, we have $(z-d)\beta + (\alpha - 1)\delta + \lambda(p-1)/d \geq (p-1)/d > l$. Equation (4) reads $\binom{(\alpha - 1)d}{(z-d)\beta + (\alpha - 1)\delta)} = 0$. But, since $(z-d)\beta < 0$, then $(z-d)\beta + (\alpha - 1)\delta < (\alpha - 1)d$. Hence $\binom{(\alpha - 1)d}{(z-d)\beta + (\alpha - 1)\delta)} \neq 0$, and the proof of Theorem 7 is complete.

\end{enumerate}

\begin{cor}
Let $f(x)$ be  a permutation binomial over $\mathbb{F}_{p}$. Then $p-1 \leq d ( d -2)$ except possibly in the case $d \equiv 0 \pmod{3}$, $p =  d^{2} - d + 1$ and one of the two possibilities: $n \equiv 0 \pmod{(p-1)/d}$ or $m \equiv 0 \pmod{(p-1)/d}$.
\end{cor}

\begin{proof}
Since there is no permutation binomial over $\mathbb{F}_2$ and over $\mathbb{F}_3$, we may suppose that $p\geq 5$. From Theorem 7, we have $(p-1)/d \leq d-1$. It remains to consider the case $(p-1)/d = d-1$, i.e. $p = d^{2}-d+1$. Suppose that there exists a permutation binomial over $\mathbb{F}_{p}$, $f(x) = ax^{n}+x^{m}$ such that  $p = d^2 - d + 1$. Consider Equation (4) for $l = d$ and let $j_{0}$ be the  integer appearinig in this equation. Since $j_{0} \in \lbrace 0, \ldots, (p-1)/d \rbrace$, then $j_{0} < d$. For any positive integer $\lambda$,  we have $j_{0} + \lambda(p-1)/d \geq j_{0} + (p-1)/d > d$ except if $j_{0} = 0$ or $j_{0} = 1$. Beyond  these exceptions,  Equation (4) reads $\binom{d}{j_{0}} = 0$. Since $j_{0} < d < p$ this equation is impossible and we get a contradiction.

\noindent (Case $j_{0} = 0$) 

Equation (4) reads $$\binom{d}{0} + \binom{d}{\frac{p-1}{d}} (a)^{\frac{p-1}{d}} = 0,$$ hence $1 + d(a)^{(d-1)} \equiv 0 \pmod p$. We deduce that $d^d\equiv(-1)^d\pmod{p}$, hence $(-d)^d\equiv 1\pmod{p}$. We have $(-d)^2\equiv d-1\pmod{p}$ and $(-d)^3\equiv 1\pmod{p}$, hence the order of $-d$ in $\mathbb{F}_p$ which is a divisor of $d$ is equal to $1$ or $3$. Since $d(d-1)=p-1$, the first possibility is excluded, hence $d\equiv 0\pmod 3$.  On the other hand side, the condition $(n-m) j_{0} + dm \equiv 0 \pmod{p-1}$ implies $dm \equiv 0 \pmod{p-1}$, i.e., $m \equiv 0 \pmod {\frac{p-1}{d}}$.

\noindent (Case $j_{0} = 1$) 

Equation (4) reads $$\binom{d}{1} + \binom{d}{1 + \frac{p-1}{d}} (a)^{\frac{p-1}{d}} = 0,$$ hence $d + (a)^{d} \equiv 0 \pmod p$. As in the preceding case we find $d^d\equiv(-1)^d\pmod{p}$, hence $d\equiv 0\pmod 3$. On the other hand side, the condition $(n-m) j_{0} + dm \equiv 0 \pmod{p-1}$ implies $(n-m) + dm \equiv 0 \pmod{p-1}$, i.e., $n \equiv 0 \pmod{\frac{p-1}{d}}$.

\end{proof}
The condition $d \equiv 0 \pmod{3}$, $p =  d^{2} - d + 1$ in Corollary 4 occurs, for instance, for $p = 7$ and $d =3$ or $p = 31$ and $d =6$ (see \cite[Cor. 2. 5]{Mas1}). This shows that the bound of $p$ in term of $d$ in Theorem 7 is sharp.
 
 If $f(x) = ax^{n} + x$ permutes $\mathbb{F}_{p}$, where $n>1$,  $a \in {\mathbb{F}_{p}}^{*}$ and $(n-1) \mid (p-1)$, then Theorem 6 does not generalize Theorem 5 which implies that $p-1 \leq (n-1)^{2}$. Theorem 7 proved above generalizes Theorem 5 even in this case as shown by the following corollary.

\begin{cor}
If $f(x) = ax^{n} + x$ permutes $\mathbb{F}_{p}$, where $n>1$ and $a \in {\mathbb{F}_{p}}^{*}$, then $p-1 \leq (n-1)(n-3)$.
\end{cor}

\begin{proof}
From Corollary 4, we have $p-1 \leq d(d-2)$, which implies that $p-1 \leq (n-1)(n-3)$, except if $d \equiv 0 \pmod {3}$, $\frac{p-1}{d} = d-1$, and $n \equiv 0 \pmod{ \frac{p-1}{d}}$ (because $m=1$). So we have only to consider the exceptionnel case. In this case, we have  $n \equiv 0 \pmod{d-1}$ and $n-1 \equiv 0 \pmod d$. Clearly $n \neq d-1$. It follows that  $n \geq 2(d-1)$. We conclude that $3 \leq d \leq \frac{n}{2} +1$. It is now easy to deduce the inequality $p-1 \leq (n-1)(n-3)$.
\end{proof}
The following result is similar to Corollary 2.4 of \cite{Mas1} except that the four polynomials arizing for $d=3$ and $p=7$ were forgoten.
\begin{cor}
If $f(x) = x^{n} + ax^{m}$ permutes $\mathbb{F}_{p}$, where $1\leq m<n<p$ and $a \in {\mathbb{F}_{p}}^{*}$, then $\gcd(n-m,p-1) > 4$  except if $d=3$, $p=7$ and $f(x)$ is one of the followings.

\noindent (i)  $f(x) = x^{4} + 3 x$.

\noindent (ii) $f(x) = x^{4} - 3 x$.

\noindent (iii) $f(x) = x^{5} +  2x^{2}$.

\noindent (iv) $f(x) = x^{5} -  2x^{2}$.
\end{cor}
\begin{proof}
We conclude from Corollary 4, that if $d = 4$, then $p-1 \leq 8$, i.e., $p \leq 7$. We see from table 7.1 of \cite{Lid2} that there are no permutation binomials in this case. When $d =2$, we conclude from Corollary 4, that there are no permutation binomials in this case. When $d =3$, Corollary 4 implies that $p =7$. We see from table 7.1 of \cite{Lid2} that the only possible cases are the one listed above.
\end{proof}
\section{Permutation binomials over a subfield of $\mathbb{F}_{q}$ arising from permutation binomials over $\mathbb{F}_{q}$}

Before stating a result about the possibilty to deduce, in some cases, a permutation binomial of a subfield of $\mathbb{F}_q$ from a given permutation binomial of $\mathbb{F}_{q}$, we make the following definition.
\begin{definition}Fix the integers $m$ and $n$ such that $1\leq m<n\leq q-1$ and let $d=gcd(n-m,q-1)$. We say that the polynomials $f(x)=ax^n+x^m$ and $g(x)=bx^n+x^m$, with coefficients in $\mathbb{F}_{q}$, are $d$-equivalent and we write $f\stackrel{d}{\sim}g$ if and only if there exists $\epsilon\in(\mathbb{F}_{q})^d$ such that $b=\epsilon a$.
\end{definition}
Obviously the above relation in the set of binomials over $\mathbb{F}_{q}$, of degree at most $q-1$, where the couple $(m,n)$ is fixed, is an equivalence relation and that each equivalence class contains $(q-1)/d$ elements.
\begin{lemma}Suppose that the polynomials $f(x)=ax^n+x^m$ and $g(x)=bx^n+x^m\in\mathbb{F}_{q}[x]$, are $d$-equivalent and that $f(x)$ permutes $\mathbb{F}_{q}$, then so does $g(x)$.
\end{lemma}
\begin{proof} Since $gcd(n-m,q-1)=d$, there exit two integers $u$ and $v$ such that
\begin{equation}
u(n-m)+v(q-1)=d.
\label{eq:9}\end{equation}
The binomials $f(x)$ and $g(x)$ being $d$-equivalent, there exists $\eta\in\mathbb{F}_{q}$ such that $b=\eta^da$. Using $(9)$, we obtain $b=\eta^{u(n-m)}a$. We deduce that 
\begin{align*}
g(x)&=\eta^{u(n-m)}ax^n+x^m=\eta^{-um)}[\eta^{un}ax^n+\eta^{um}x^m]\\
&=\eta^{-um}[a(\eta^{u}x)^{n}+(\eta^{u}x)^m]=\eta^{-um}f(\eta^{u}x)),
\end{align*}
and this proves our lemma.
\end{proof}

\begin{thm}
\label{theorem 3}
Let $f(x) = ax^{n} + x^{m}$ be a permutation binomial of $\mathbb{F}_{q}$ with $q=p^r$ and $s$ be a positive divisor of $r$. Let $d = {\mbox{gcd}}(n-m,q-1)$. 
\begin{itemize}
\item[( 1)] There exists a binomial $g(x) =b x^{n} + x^{m}\in\mathbb{F}_{p^s}[x]$ $d$-equivalent to $f(x)$ if and only if the order of $a$ in $(\mathbb{F}_{q})^{\star }$ divides $lcm(p^s-1, (q-1)/d).$
\item [(2)]If these equivalent conditions hold, then the number of $g(x) =b x^{n} + x^{m}\in\mathbb{F}_{p^s}[x]$, $d$-equivalent to $f(x)$ is equal to $gcd(p^s-1,(q-1)/d)$ and they are all distinct as permutations of  $\mathbb{F}_{p^s}$. Moreover, we have $g(x) \equiv bx^{n_{1}} + x^{m_{1}} \pmod{x^{p^s} - x}$ if 
$p^s - 1 \nmid d$ and $g(x) \equiv (b+1)x^{k} \pmod {x^{p^s}-x}$ if $p^s-1 \mid d$ where $k, m_{1}, n_{1}$ are positive integers less than $p^s-1$, $m_{1} \neq n_{1}$, $\gcd(p^s-1, k) = 1$.
\item[(3)]Let $t$ be a posive integer (not necesseraly dividing $r$). There exists a binomial $g(x) =b x^{n} + x^{m}\in\big(\mathbb{F}_{p^t}\cap\mathbb{F}_{q}\big)[x]$, $d$-equivalent to $f(x)$ if and only if the order of $a$ in $(\mathbb{F}_{q})^{\star }$ divides $lcm(p^t-1, (q-1)/d).$
\end {itemize}
\end{thm}
\begin{proof}
\noindent (1) Suppose first that the order of $a$ in $(\mathbb{F}_{q})^{\star }$ divides $lcm(p^s-1, (q-1)/d).$ We will use the following claim for which the proof is omited

\textbf{Claim 1} Let $\delta, u, v$ be positive integers. Then $\delta \mid lcm(u, v)$
if and only if there exist positive integers $\delta_{1}, \delta_{2}$ such that
$\delta_{1} \mid u$, $\delta_{2} \mid v$ and $\delta = lcm(\delta_{1}, \delta_{2})$.

Let $\delta$ be the order of $a$ in $\mathbb{F}_{q}^{\star}$, then $\delta = lcm(\delta_{1}, \delta_{2})$, where
$\delta_{1}$ and $\delta_{2}$ are positive integers such that $\delta_{1} \mid p^s-1$
and $\delta_{2} \mid (q-1)/d$. Let $\xi$ be a generator of $\mathbb{F}_{q}^{\star}$, then $a
= (\xi^{(q-1)/\delta_{1}})^{i}(\xi^{(q-1)/\delta_{2}})^{j}$ for some nonnegative
integers $i$ and $j$. Let $\epsilon = \xi^{-j(q-1)/\delta_{2}}$, $b = \epsilon
a$ and $g(x)=b x^{n} + x^{m}$. Then $\epsilon^{(q-1)/d} = (\xi^{-j(q-1)/d\delta_{2}})^{q-1} = 1$ and $b^{p^s-1} = (\xi^{i(p^s-1)/ \delta_{1}})^{q-1} = 1$, hence $\epsilon\in\mathbb{F}_{q})^d$, $b\in\mathbb{F}_{p^s}$ an $g(x)\stackrel{d}{\sim}f(x)$.

Conversely, suppose that there exist $g(x)=b x^{n} + x^{m}\in\mathbb{F}_{p^s}$, $d$-equivalent to $f(x)$, then we may find a $(q-1)/d$-th root of unity $\epsilon$ such that $a=\epsilon b$. Hence $a^{ lcm(p^{s}-1, (q-1)/d)}=(\epsilon b) ^{\frac{(p^s-1)(q-1)/d}{\delta_s}}=1$, where $\delta_s=gcd(p^s-1,(q-1)/d)$. It follows that the order of $a$ in $(\mathbb{F}_{q})^{\star }$ divides $lcm(p^s-1, (q-1)/d)$.

\noindent (2) Let $\delta_s=gcd(p^s-1, (q-1)/d)$. By (1), there exists at least one permutation binomial of $\mathbb{F}_{q}$, $g(x)=c_s x^{n} + x^{m}$ with $c_s\in\mathbb{F}_{p^s}$, $d$-equivalent to $f(x)$. Let $h(x)=b_s x^n+x^m$ be any permutation binomial of $\mathbb{F}_{q}$with $b_s\in\mathbb{F}_{p^s}$, $d$-equivalent to $f(x)$, then $g(x)\stackrel{d}{\sim}h(x)$, hence there exists $\epsilon\in\mathbb{F}_{q})^d$ such that $b_s=\epsilon c_s$. We deduce that $\epsilon=b_s/c_s\in\mathbb{F}_{p^s}$. It follows that $\epsilon^{p^s-1}=1=\epsilon^{(q-1)/d}$ and then $\epsilon^{\delta_s}=1$. We conclude that $h(x)$ has the form $h(x)=\epsilon c_s x^n+x^m$ with $\epsilon$ satsfying the condition $\epsilon^{\delta_s}=1$. On the other hand any polynomial $h(x)$ of this form is $d$-equivalent to $g(x)$ and then to $f(x)$. Clearly all these $h(x)$, as permutations of $\mathbb{F}_{q}$, are distinct. Because all of them take different values at the argument $x=1$, they are distinct as permutations of $\mathbb{F}_{p^s}$. We conclude that the number of these $h$'s is equal to $\delta_s$.

To prove the last part of the theorem we reduce $g(x)$ modulo $x^{p^s} - x$. Denote by $\overline{g(x)}$
the unique polynomial over $\mathbb{F}_{p^s}$ of degree at most $p^s-1$ such that $g(x)\equiv\overline{g(x)}\pmod{x^{p^s} - x}$. Set $n =
(p^s-1)\lambda + n_{1}$ and $m = (p^s-1)\mu + m_{1}$ with $0 \leq m_{1}, n_{1}\leq p^s-2$.
If $m_{1} = 0$ or $n_{1} = 0$, then the degree of $\overline{g(x)}$ is equal to
$p^s-1$, which is excluded by the fact that $g(x)$ is a permutation polynomial of
$\mathbb{F}_{p^s}$. If $m_{1} = n_{1}$, then clearly $p^s - 1 \mid d$ and
$\overline{g(x)} = (b+1)x^{k}$, where $k = n_{1} = m_{1}$ and $\gcd(k, p^s-1) = 1$.
Suppose now that $m_{1} \neq 0$, $n_{1} \neq 0$ and $m_{1} \neq n_{1}$, then $p^s-1
\nmid d$ and  $\overline{g(x)} = bx^{n_{1}} + x^{m_{1}}$. Let $k = \gcd(m_{1},
n_{1})$, then the polynomial $g_{1}(x) = bx^{n_{1}/k} + x^{m_{1}/k}$ is a permutation
binomial of $\mathbb{F}_{p^s}$.
\noindent (3) We will use the following which is certainly well known.

\textbf{Claim 2} Let $a$, $b$, $c$ be nonzero integers, then $$gcd\big(lcm(a,b),lcm(a,c)\big)=lcm\big(a,gcd(b,c)\big)$$.

Suppose that there exists a binomial $g(x) =b x^{n} + x^{m}\in\big(\mathbb{F}_{p^t}\cap\mathbb{F}_{q}\big)[x]$, $d$-equivalent to $f(x)$. Let $s=gcd(r,t)$, then $\mathbb{F}_{q}\cap\mathbb{F}_{p^t}=\mathbb{F}_{p^s}$ and by $(1)$, the order of $a$ in $\mathbb{F}_{q}$ divides $lcm(p^s-1, (q-1)/d)$. Therefore this order divides $lcm(p^t-1, (q-1)/d)$. Suppose now that the order of $a$ in $\mathbb{F}_{q}$ divides $lcm(p^t-1, (q-1)/d)$, then applying the above Claim with $a=(q-1)/d$, $b=p^t-1$ and $c=p^r-1$, we conclude that this order divides $lcm(p^s-1, (q-1)/d)$ and then by $(1)$ there exists a binomial $g(x) =b x^{n} + x^{m}\in\big(\mathbb{F}_{p^t}\cap\mathbb{F}_{q}\big)[x]$, $d$-equivalent to $f(x)$.

\end{proof}
\begin{rem}
Suppose that $p$ is odd, then under the hypothesis $(1)$ of the above theorem, we have
$\gcd(d, p^s-1) \neq 1$. Indeed if this gcd is equal to 1, then $lcm(p^s-1, (q-1)/d) =
(q-1)/d$. But it is knownn that if $(-1/a)^{(q-1)/d} = 1$, then the corresponding
binomial is not a permutation binomial of $\mathbb{F}_{q}$(see \cite{Sma2}).
\end{rem}
\begin{corollary}
Let $f(x) = ax^{n} + x^{m}$ be a permutation binomial of $\mathbb{F}_{q}$. Suppose
that the order of $a$ in $\mathbb{F}_{q}^{*}$ divides $lcm(p-1, (q-1)/d).$ Then $p-1 \leq d(d-1)$.
\end{corollary}

\begin{proof} If $p-1\mid d$, the corollary is clear. If not, the proof is a direct consequence of Theorem 7 and Theorem 9.
\end{proof}
\begin{corollary} Suppose that there exists a permutation binomial $f(x)=ax^n+x^m$ of $\mathbb{F}_q$ with q= $p^r$ such that for any prime number, $l\mid d$, we have $\mbox{gcd}(l(l-1),r)=1$. Then $d=p-1$ or there exist a permutation binomial of $\mathbb{F}_p$, $g_1(x)=cx^{n_1}+x^{m_1}$ such that $n\equiv kn_1\pmod{p-1}$, $m\equiv km_1\pmod{p-1}$, $0<km_1<kn_1<p-1$, where $k$ is a positive integer coprime with $p-1$, and $p-1 \leq d(d-1)$. Moreover the two possibilites exclude each other.
\end{corollary}

\begin{proof} Let $l$ be any prime factor of $d$. We have $p^r\equiv 1\pmod {l}$, by asymption and $p^{l-1}\equiv 1\pmod {l}$, by Fermat's little theorem. Since $r$ and $l-1$ are coprime, then $p\equiv 1\pmod l$. It is easy to see that $p\equiv 1\pmod d$ and $\mbox{lcm}(p-1,(q-1)/d)=q-1$ so that Theorem 9 may be applied to any permutation binomial of $\mathbb{F}_q$. Let $g(x)$ be the permutation binomial of $\mathbb{F}_q$ with coefficients in  $\mathbb{F}_p$ deduced from $f(x)$, using Theorem 9. Let $g_1(x)$ be the reduced polynomial of $g(x)$ modulo $x^p-x$. Then $g_1(x)$ is a monomial or $g_1(x)$ is a sum of $2$ monomials of degree $n'$ and $m'$ respectively satisfying $0<m'<n'<p-1$. Moreover the first case holds if and only if $p-1\mid d$. Since $p\equiv 1\pmod d$ then the first case holds if and only if $d=p-1$. To complete the proof let $k=\mbox{gcd}(m',n')$, $m_1=m'/k$ and $n_1=n'/k$, and by applying Corollary 6, we have $p-1 \leq d(d-1)$.
\end{proof}
The following result is a generalization of \cite[Corallary 2. 4, Corollary 2. 5]{Mas1}.
\begin{corollary} There does not exist a permutation binomial of $\mathbb{F}_q$ with $q=p^r$ if one of the following conditions holds.
\begin{itemize}
\item[(i)] $r$ odd, $d=2$, $p\neq3$.
\item[(ii)] $r$ odd, $d=4$, $p\neq5$.
\item[(iii)] $\mbox{gcd}(r,6)=1$, $d=3$, $p\neq 7$.
\item[(iv)] $\mbox{gcd}(r,10)=1$, $d=5$ $p\neq 11$.
\item[(v)] $\mbox{gcd}(r,6)=1$, $d=6$, $p\neq 7, 13, 19, 31$.
\item[(vi)] $\mbox{gcd}(r,42)=1$, $d=7$, $p\neq 29$.
\item[(vii)] $r$ odd, $d=8$, $p\neq 17$.
\end{itemize}
\end{corollary}
\begin{proof} We prove the case $(v)$ using  corollaries 6 and 8 and \cite[Corollary 2. 5]{Mas1}. The proof of the other statements will be omitted. Suppose that there exist a permutation binomial of $\mathbb{F}_q$ with $d=6$ and $\mbox{gcd}(r,6)=1$, then the hypotheses of the above corollary holds. We deduce that $p=d+1=7$ or there exists some permutation binomial $g_1(x)=cx^{n_1}+x^{m_1}$ of $\mathbb{F}_p$. It is evident that $\mbox{gcd}(n_1-m_1,p-1)$ divides $d=6$ and is not trivial. The possible values of this gcd are $2$ or $3$ or $6$. According to \cite[ Corollary 2. 5]{Mas1}, the possible values of $p$ are $p=7, 13, 19, 31$.
\end{proof}
\begin{rem} It is of interest to improve the conditions on $r$ and $p$ in the above corollary. We use the results of \cite[Table 7.1] {Lid2} to 
make some observations in this direction. Since $ax^3+x$ is a permutation polynomial of $\mathbb{F}_q$ for $q\equiv 0\pmod{3}$ and $-a$ not a square, then the condition $p\neq3$ is necessary for $d=2$. The polynomial $ax^5+x$ is a permutation of $\mathbb{F}_q$ for $q\equiv 0\pmod{5}$ and $-a$ is not fourth power, hence the condition $p\neq5$ is necessary for $d=4$. Let $a\in\mathbb{F}_9$ such that $a^2=-1$, then $ax^5+x$ permutes $\mathbb{F}_9$, hence the condition $r$ odd is necessary for $d=4$.
\end{rem}
\begin{proposition} Le $q=p^r$ and $f(x) =ax^{n} + x^{m}\in\mathbb{F}_q[x]$ be a permutation binomial. Let $\mathbb{F}_{p^{s_1}},\ldots,\mathbb{F}_{p^{s_u}}$ be subfields of $\mathbb{F}_q$ such that for each $i$, $\mathbb{F}_{p^{s_i}}$ contains the coefficients of some binomial $g_i(x)$, $d$-equivalent to $f(x)$, then $\cap_{i=1}^{u}\mathbb{F}_{p^{s_i}}$ contains the coefficients of some binomial $g(x)$, $d$-equivalent to $f(x)$.

\end{proposition}
\begin {proof} We prove the result for $u=2$. The proposition may be completed easily by induction. By Theorem 9, $(1)$, the order of $a $ divides both  $lcm(p^{s_1}-1,(q-1)/d)$ and $lcm(p^{s_2}-1,(q-1)/d)$, hence by Claim 2, this order divides $lcm\big((q-1)/d,gcd(p^{s_1}-1,p^{s_2}-1)\big)$. It is well known that $gcd(p^{s_1}-1,p^{s_2}-1)=p^{gcd(s_1,s_2)}-1$ and that $\mathbb{F}_{p^{s_1}}\cap\mathbb{F}_{p^{s_2}}=\mathbb{F}_{p^{gcd(s_1,s_2)}}$. By Theorem 9 again this last field contains the coefficients of some binomial $g(x)$, $d$-equivalent to $f(x)$.

\end{proof}

If we consider all the subfields $\mathbb{F}_{p^{s_i}}$ of $\mathbb{F}_q$ satisfying the given property in the preceding proposition we may conclude that the field $F_0=\cap_{i}\mathbb{F}_{p^{s_i}}$ contains the coefficients of some binomial $g(x)$, $d$-equivalent to $f(x)$. We call this field {\itshape the smallest field containing the coefficients of some $d$-equivalent to $f(x)$.}

The next proposition shows that binomials that are conjugate over $\mathbb{F}_q$ or in the same $d$-class have the same smallest field.
\begin{proposition} Let  $f(x) =ax^{n} + x^{m}\in\mathbb{F}_q[x]$ be a permutation binomial of $\mathbb{F}_q$ and let $F_0$ be the smallest field containing the coefficients of some $d$-equivalent to $f(x)$.
\begin{itemize}
\item[(1)] Let $g(x)\in\mathbb{F}_q[x].$ If $f\stackrel{d}{\sim}g$, then $F_0$ is the smallest field containing the coefficients of some $d$-equivalent to $g(x)$.
\item[(2)] Let $\tilde{f}(x) =a^{p^e}x^{n} + x^{m}$, then $F_0$ is the smallest field containing the coefficients of some $d$-equivalent to $\tilde{f}(x)$.
\end{itemize}
\end{proposition}
\begin{proof}  $(1)$ Let $F_1$ be the smallest field corresponding to $g(x)$. For the proof of $(1)$ and by symmetry it is sufficient to prove that $F_0\subset F_1$. Let $g_1(x)$ be a $d$-equivalent of $g(x)$ with coefficients in $F_1$, the $g_1\stackrel{d}{\sim}g\stackrel{d}{\sim}f$, hence $F_1$ contains the coefficients of some $d$-equivalent to $f(x)$. Therefore $F_0\subset F_1$.

$(2)$ The result follows from Theorem 9, $(1)$ and the observation that $a$ and $a^{p^e}$ have the same order.

\end{proof}

\end{document}